\documentclass[conference]{IEEEtran}
\IEEEoverridecommandlockouts
\usepackage{cite}
\usepackage{amsmath,amssymb,amsfonts}
\usepackage{algorithmic}
\usepackage{algorithm}
\usepackage{array}
\usepackage{graphicx}
\usepackage{textcomp}
\usepackage{xcolor,ulem}
\usepackage{rotating}
\usepackage{hyperref}
\usepackage{leftidx}
\usepackage{multirow}
\usepackage{comment}
\usepackage[T1]{fontenc}
  \usepackage{ulem}
\def\BibTeX{{\rm B\kern-.05em{\sc i\kern-.025em b}\kern-.08em
    T\kern-.1667em\lower.7ex\hbox{E}\kern-.125emX}}

\newtheorem{theorem}{Theorem}

\newcommand{\yhk}{z_{k}}
\newcommand{\yhkun}{z_{k+1}}

\newcommand{\xhk}{x_{k}}
\newcommand{\xhkun}{x_{k+1}}

\newcommand{\R}{R}

\newcommand{\vHk}{v_{H}}

\newcommand{\RR}{\mathbb{R}}

\newcommand{\Id}{\mathrm{Id}}

\DeclareMathOperator*{\argmin}{arg\,min}
\newcommand{\minimize}[1]{\ensuremath{\underset{\substack{{#1}}}{\mathrm{minimize}}\;\; }}
\DeclareMathOperator*{\Argmin}{Argmin}
\newcommand{\prox}{\normalfont \textrm{prox}}

\newcommand{\ftn}{\scriptsize}

\begin{document}

\title{%
A multilevel framework for accelerating \\
uSARA in radio-interferometric imaging
\thanks{This work is funded the Fondation Simone et Cino Del Duca - Institut de France. The authors thank the Centre Blaise Pascal of ENS Lyon for the computation facilities. The platform uses SIDUS \cite{quemener2013}, which was developed by Emmanuel Quemener.}
}

\author{\IEEEauthorblockN{Guillaume Lauga$^\dagger$, Audrey Repetti$^\ddagger$$^{**}$, Elisa Riccietti$^\dagger$, Nelly Pustelnik$^{*}$, Paulo Gonçalves$^\dagger$,  Yves Wiaux$^{**}$}
\IEEEauthorblockA{$^\dagger$\textit{Univ Lyon, Inria, EnsL, UCBL, CNRS, LIP, UMR 5668, F-69342} \\
$^{*}$\textit{Ens de Lyon,
		CNRS, Laboratoire de Physique, F-69342}\\
$^\ddagger$\textit{Department of Actuarial Mathematics \& Statistics, Heriot-Watt University, Edinburgh EH14 4AS, United Kingdom} \\
$^{**}$ \textit{Institute of Sensors, Signals and Systems, Heriot-Watt University, Edinburgh EH14 4AS, United Kingdom} \\
\{guillaume.lauga,elisa.riccietti,nelly.pustelnik,paulo.goncalves\}@ens-lyon.fr ; \{a.repetti,y.wiaux\}@hw.ac.uk}
}
\maketitle

\begin{abstract}
This paper presents a multilevel algorithm specifically designed for radio-interferometric imaging in astronomy. The proposed algorithm is used to %
solve the uSARA (unconstrained Sparsity Averaging Reweighting Analysis \cite{carrillo2012sparsity, repetti2020forward,terris2023RI}) formulation of this image restoration problem. Multilevel algorithms rely on a hierarchy of approximations of the objective function to accelerate its optimization. In contrast to the usual multilevel approaches where this hierarchy is derived in the parameter space, here we construct the hierarchy of approximations in the  observation space. The proposed approach is compared to a reweighted forward-backward procedure, which is the backbone iteration scheme for solving the uSARA problem. %
\end{abstract}


\section{Introduction}
Radio-interferometric (RI) imaging aims to reconstruct a sky brightness distribution from noisy Fourier observations (named visibilities). The Fourier frequency sample distribution is dictated by the position of the antennas used to probe the sky. The number of antennas within an interferometer being finite, this leads to undersampling the Fourier measurements \cite{thompson2017interferometry} (see Figure \ref{fig:uv_coverage}). 
Obtaining the target image requires to %
deconvolve and denoise the Fourier measurements.%

The goal is thus to  reconstruct the intensity image $\overline{x} \in \RR^n$ from a set of measured complex \textit{visibilities} $y\in \mathbb{C}^m$ acquired in the Fourier space. The corresponding discretized forward model can be formulated as follows \cite{carrillo2012sparsity}: 
\vspace{-0.2em}
\begin{equation} \label{pb:inv}
    y = \mathbf{GFZ}\overline{x} + \epsilon:= \mathbf{\Phi}\overline{x} + \epsilon\vspace{-0.2em}
\end{equation} 
where $\mathbf{G}\in \mathbb{C}^{m\times d}$ is a sparse interpolation operator %
that maps the Fourier coefficients on a regular grid to the non-uniformly located visibilities%
, $\mathbf{F} \in \mathbb{C}^{d \times d}$ is the $2$D Discrete Fourier Transform, and $\mathbf{Z} \in \RR^{d \times n}$ is a zero-padding operator to properly map $\overline{x}$ for the convolution performed through the operator $\mathbf{G}$. %
The term $\epsilon$ stands for a centered white Gaussian noise\footnote{This modelling of the measurement operation is an approximation of the true measurement process \cite{onose2016scalable,mhiri2023contributions,birdi2019advanced}: we assume here a narrow field of view, white noise across all visibilities and no anisotropic perturbations.}. Typically %
$m\gg n$ and is of the order of 10 millions of visibilities (see \cite{mhiri2023contributions} for more details).%

This problem presents two main challenges. First, it is severely ill-posed, hence requiring advanced imaging techniques. Second, the size of the data streams coming from radio-telescopes is expected to be ever increasing, thus raising the challenge of designing highly scalable methods. %

\noindent \textbf{Recovery techniques in RI} -- 
CLEAN \cite{hogbom1974aperture,cornwell2008multiscale} is the most used algorithm in RI imaging. %
It is similar to a matching pursuit method, %
but shows limitations when probing extended complex emissions%
, or for large numbers of point sources %
\cite{yatawatta2010fundamental}. %
Penalized variational procedures have hence been proposed to improve the quality of the reconstructed images in a more general framework \cite{mhiri2023contributions,birdi2019advanced}.
The state-of-the-art variational formulation in RI is SARA \cite{carrillo2012sparsity}%
, promoting average sparsity of the solution in a concatenation of bases through a reweighted-$\ell_1$ procedure. %
Specifically, it aims at minimizing 
a log-sum prior by solving a sequence of weighted $\ell_1$ problems  \cite{candes2008,repetti2021variable,repetti2020forward}. 
Two SARA formulations have been proposed for imaging \cite{carrillo2012sparsity, repetti2020forward, terris2023RI}: a constrained and an unconstrained one aka uSARA, on which we will focus in this work. 

\noindent \textbf{Scaling to high-dimensional data} -- 
Many algorithms have been proposed to reduce the computational load induced by the number of visibilities. Most often, these algorithms consider only a subset of visibilities 
at each iteration. 
A first idea is to split the visibilities into blocks and to parallelize the action of $\mathbf{\Phi}$ blockwise \cite{onose2016scalable}. 
In \cite{thouvenin2023parallel}, this approach was also extended to multi-spectral data  with a modification of SARA to take into account spectral correlations. 
Other approaches solve approximations of the original problem in smaller dimensions, by selecting relevant visibilities in a sketching fashion (suitable random projections for instance) \cite{vijay2017fourier} %
or updating only with a fraction of the visibilities at each optimization step in an online manner \cite{cai2019online}. %
Finally acceleration schemes such as preconditioning strategies can be considered to better take into account the specific RI Fourier distribution \cite{onose2017accelerated}.

\noindent \textbf{Contributions} -- 
We leverage a recent work \cite{lauga2023} of the authors to design an accelerated algorithm that combines a \textit{multilevel} (ML) procedure with FISTA iterations \cite{beck2009, chambolle2015}
for solving the uSARA problem. 
%
ML algorithms have been shown to significantly accelerate the resolution of problems whose structure endow a hierarchy of functions defined on nested subspaces of the parameters space, able to approximate the objective function \cite{parpas2016,parpas2017,ang2023,lauga2022,lauga2023}.
%
%
%
%
In this work, we propose an original ML framework adapted to RI, where lower dimensional nested subspaces are constructed in the visibility domain, rather than in the image domain as usually done for ML approaches.
We show through simulations that the resulting ML  %
scheme yields 
significant acceleration when solving the uSARA problem. %

\section{The multilevel framework}
\label{sec:multilevel_framework}

Inverse problems of the form of~\eqref{pb:inv} can be solved by defining iterations to \vspace{-0.7em}
\begin{equation}\label{min:fine-gen}
    \minimize{x \in \mathbb R^n} F(x):= \underset{L(x)}{\underbrace{\frac12 \| \mathbf{\Phi}x -y \|^2 }} + R(x), \vspace{-0.5em}
\end{equation}
where $R \colon \mathbb R^n \to ]-\infty, +\infty]$ is a regularization function incorporating prior information on the target solution.

Without loss of generality, we present the proposed ML strategy on a two-level case.
In this setting, 
we index %
the functions at the coarse level with  subscript $H$, i.e. $F_H$, $L_H$ and $R_H$, for $F$, $L$, and $R$, respectively. \smallskip

\noindent \textbf{Spirit of the method} -- Given an objective function at the fine level $F$, the goal of ML approaches is to build a coarse approximation $F_H$ of $F$, which is cheaper to optimize, to accelerate the minimization of $F$. 
Then, a ML algorithm consists of alternating iterations at the coarse level on $F_H$ ($\mathbf{ML}$ steps) and at the fine level on $F$. Within the ML framework developed in~\cite{lauga2023}, the minimization of~\eqref{min:fine-gen} at the fine level can be computed using either forward-backward iterations or its accelerated inertial version FISTA. %
Then, the overall ML alternating procedure reads
\begin{equation}
\begin{array}{l}
\text{for } k=0, 1, \ldots\\
\left\lfloor 
\begin{array}{l}
    \overline{x}_{k} = \mathbf{ML}(x_{k}) \\
    \xhkun = \mathbf{FISTA}(\overline{x}_{k}).
\end{array}
\right.
\end{array}
\end{equation}
The %
crucial component of a ML strategy is the construction of a $F_H$ in order to be consistent with $F$ %
\cite{parpas2017,lauga2023}. 

\smallskip
The main contribution of this article is to construct $F_H$ exploiting properties of the RI problem, considering a coarse model in the data domain, and leveraging the specific RI Fourier sub-sampling.

\subsection{Proposed coarse model in data space}
The usual approach to construct $F_H$ consists in approximating $F$ in a lower dimensional space.
This would amount here to decrease the size of the image $x$ and to formulate a similar optimization problem for a low resolution image \cite{parpas2017,lauga2023}. 
However, to take into account that the limiting factor in RI imaging is the large number of visibilities rather than the size of the sought image, we deviate from the classical ML scheme, and we construct a coarse model based on the following approximation of $L$:
\begin{equation}
    \label{eq:sketching_model}
   (\forall x \in \RR^n) \quad L_H(x) :=\frac{1}{2}\Vert \mathbf{S}\mathbf{\Phi} x-\mathbf{S}y \Vert^2 
\end{equation}
where $\mathbf{S} \colon \mathbb C^m \to \mathbb C^{m_H}$ is an operator reducing the data dimensionality $m$ to a lower dimension $m_H < m$. %
Note that $L_H \colon \mathbb R^n \to \mathbb R$ is defined on the same space $\mathbb R^n$ as $L$, thus information transfer operators between levels, commonly used in standard multilevel algorithms \cite{parpas2016,parpas2017,ang2023,lauga2022,lauga2023},  are not required in the proposed setting.%

Formulation \eqref{eq:sketching_model} is standard in the sketching literature, where $\mathbf{S}$ is typically a Gaussian random matrix so that minimizing \eqref{eq:sketching_model} guarantees signal recovery \cite{foucart2013}. Such a strategy has however many drawbacks for RI imaging 
that were investigated in \cite{vijay2017fourier}. Notably operator $\mathbf{S}\Phi$ is dense and thus computationally intensive to apply in iterative optimization.
In general, choosing $\mathbf{S}$ to reduce computation complexity without sacrificing reconstruction accuracy is challenging, and some choices may lead to sub-optimal reconstruction \cite{vijay2017fourier}. However %
in our framework, $L_H$ is only used to propel the minimization of the fine level objective function. %
We propose to sub-sample the $u$-$v$ coverage, that will enable preserving the reconstruction quality while reducing the computation complexity of the overall minimization method. This choice will be further discussed in Section~\ref{sec:alg}. \vspace{-0.5em}%

\section{ML approach for uSARA acceleration}

\subsection{uSARA approach in a nutshell}

The uSARA problem can be written as in \eqref{min:fine-gen}, where $R$ corresponds to a log-sum penalization to promote sparsity in the concatenation of the first eight Daubechies wavelet bases and the Dirac basis. 
Such a regularization is then handled using a reweighted $\ell_1$ approach, that aims at solving a sequence of weighted $\ell_1$ problems \cite{candes2008, repetti2020forward, terris2023RI}. 
The resulting reweighting procedure can then be written as \vspace{-0.35em}
\begin{equation}
    \label{eq:reweight_step}
\begin{array}{l}
    \text{for } i = 0, \ldots, I \\
    \left\lfloor
    \begin{array}{l}
        \displaystyle \widetilde{x}_{i+1}  \in  \Argmin_{x \in \RR^n} F_i(x):= L(x) + R(x,\mathbf{W}_i), \\
        \mathbf{W}_{i+1}  = \text{Diag}\left(\lambda (\rho + |\mathbf{\Psi}^*\widetilde{x}_{i+1}|)^{-1}\right), 
    \end{array}
    \right.
\end{array}
\end{equation}
where $I>0$ is the maximum number of reweighting steps, $\R(x,\mathbf{W}_i) = \Vert \mathbf{W}_i \mathbf{\Psi}^* x \Vert_1 + \iota_{\RR^n_+}(x)$, with $\mathbf{\Psi}$ being the SARA dictionary, and $\iota_{\RR^n_+}$ being the indicator function associated with the positive orthant, 
$\lambda>0$ is a regularization parameter balancing the contribution of the regularization and the data fidelity terms, and $\rho>0$ ensures stability of the method. 
As $\rho$ tends to $0$ the solution of the weighted $\ell_1-$norm problem approaches that of the $\ell_0$ pseudo-norm problem.  

It has been shown in~\cite{repetti2021variable} that when solving approximately the minimization problem in~\eqref{eq:reweight_step} with a fixed number of forward-backward iterations, the sequence $(\widetilde{x}_i)_{i \in \mathbb N}$ converges to a critical point of $F$ in \eqref{min:fine-gen}. \vspace{-0.7em}

\subsection{Proposed IML-FISTA for uSARA}
In the context of RI imaging we adapt the inexact ML FISTA (IML-FISTA) proposed in \cite{lauga2023} for minimizing $F_i$ at each reweighting step $i\in \{0, \ldots, I\}$. %
Then, at the $i$-th reweighting step, the proposed IML-FISTA iterations for uSARA read
\begin{equation}\label{algo:full-rw}
\begin{array}{l}
\text{for } k = 0, 1, \ldots \\
\left\lfloor
\begin{array}{l}
    \overline{z}_{k} = \mathbf{ML}(\yhk) \\
    \xhkun  \approx
        \prox_{\tau R( \cdot, \mathbf{W}_i)}\left(\overline{z}_{k} - \tau \nabla L(\overline{z}_{k})\right),  \\
    \yhkun  = \xhkun + \alpha_{k}(\xhkun-\xhk),
\end{array}
\right.
\end{array}
\end{equation}
where $\tau>0$ and $(\alpha_{k})_{k \in \mathbb N}$ are chosen according to \cite{aujol2015}, and the approximation errors on the proximal operator are assumed to be summable \cite{aujol2015,lauga2023}. 
The multilevel (\textbf{ML}) step consists of updating the variable $\yhk$ at certain iterations with a correction from coarse models to obtain a better update $\overline{z}_{k}$. %
The detailed version of this step and the variables involved are presented in Algorithm~\ref{alg:IMLFISTA}. 

\begin{algorithm}
\begin{algorithmic}
  \IF{Coarse correction at iteration $k$}
   \STATE Set $\tau_H>0$ and $\alpha_{H}>0$ according to \cite{aujol2015}
\label{alg_step:projection}\\ %
    \STATE $z_+  = \underbrace{(\Id-\tau_H\nabla F_H)\circ \ldots \circ (\Id-\tau_H\nabla F_H)}_{p \text{ times}}(z_k) $  \label{alg_step:minimization}
    \STATE $\overline{z}_{k} = z_k + \alpha_{H}\left(z_+-z_k\right)\label{alg_step:prolongation}$ %
   
    \ELSE
     \STATE $\overline{z}_{k} = z_{k}$ \label{alg_step:standard_update}
     \ENDIF
     
\end{algorithmic}
  \caption{ML step in~\eqref{algo:full-rw} \label{alg:IMLFISTA} }
\end{algorithm}\vspace{-0.5em}

\noindent At iteration $k$ of algorithm~\eqref{algo:full-rw} the coarse objective function $F_H$ is given by \vspace{-0.5em}
\begin{equation}
    F_{H} = L_H + R_{H,\gamma}(\cdot,\mathbf{W}_i) + \langle \vHk, \cdot \rangle, \vspace{-0.5em}
    \label{eq:F_H}
\end{equation}
where \vspace{-0.5em}
\begin{equation}
\label{eq:v_Hk_inexact1}
\vHk = \nabla \big( L +  R_{\gamma}(\cdot, \mathbf{W}_i)
-  L_H -  R_{ H, \gamma}(\cdot, \mathbf{W}_i) \big) (\yhk).
\end{equation}
In~\eqref{eq:v_Hk_inexact1},  $R_{\gamma}(\cdot, \mathbf{W}_i)$ corresponds to a smooth approximation of $R( \cdot, \mathbf{W}_i)$ with parameter $\gamma>0$ \cite[Definition 2.1]{beck2012}, obtained using a Moreau-Yosida smoothing technique (see \cite{beck2012,lauga2023} for details). The advantage of this technique is %
that the gradient has a closed form expression. Similarly, $R_{H, \gamma}$ is a smooth approximation of the coarse approximation $R_{H}$, built using the same technique.
\noindent In~\eqref{eq:F_H}$, \vHk$ imposes first-order coherence between the smoothed versions of the functions at fine and coarse levels. %
According to \cite[Thm~2.15 and Thm~2.16]{lauga2023}, we then have the following theoretical guarantees:
\begin{theorem}
Let $i\in \{1, \ldots, I\}$. Let $(x_k)_{k\in \mathbb N}$ and $(z_k)_{k\in \mathbb N}$ be sequences generated by algorithm~\eqref{algo:full-rw}.
Assume that, for every $k\in \mathbb N$, the coarse model defined in Algorithm~\ref{alg:IMLFISTA} decreases, i.e. $F_H(z_+) \leq F_H(z_k)$\footnote{This is ensured as soon as $\tau_H < \beta_H^{-1}$, where $\beta_H>0$ is the Lipschitz constant of $\nabla F_H$.}. Then, the following assertions hold:
\begin{enumerate}
    \item $\big( F(x_k)-F^* \big)_{k \in \mathbb N}$ is decreasing at a rate of $1/k^2$,
    \item $(x_k)_{k \in \mathbb N}$ converges to a minimizer of $F_i$ when $k\to \infty$.
\end{enumerate} 
\end{theorem}

\subsection{Algorithmic settings and implementation}
\label{sec:alg}

\smallskip\noindent \textbf{Proximity operator computation} -- 
In~\eqref{algo:full-rw}, the proximity operator of $R(\cdot, \mathbf{W}_i)$ is defined, for every $x \in \RR^n$, as
\begin{equation*}
    \prox_{\tau R(\cdot, \mathbf{W}_i)}(x) = \argmin_{u \in \RR^n_+}  \frac{1}{2\tau} \Vert u-x \Vert^2 + \Vert \mathbf{W}_i \Psi^* u \Vert_1   .
\end{equation*}
Since this proximity operator does not have a closed form expression, it can be computed with sub-iterations. In particular, the dual forward-backward algorithm proposed in \cite{vu2010dualization} 
produces a sequence of feasible iterates converging to $\prox_{\tau R(\cdot, \mathbf{W}_i)}(x)$ \cite[Thm 3.7] {vu2010dualization}.
%
%

\begin{figure}
    \centering    \includegraphics[width=0.25\textwidth]{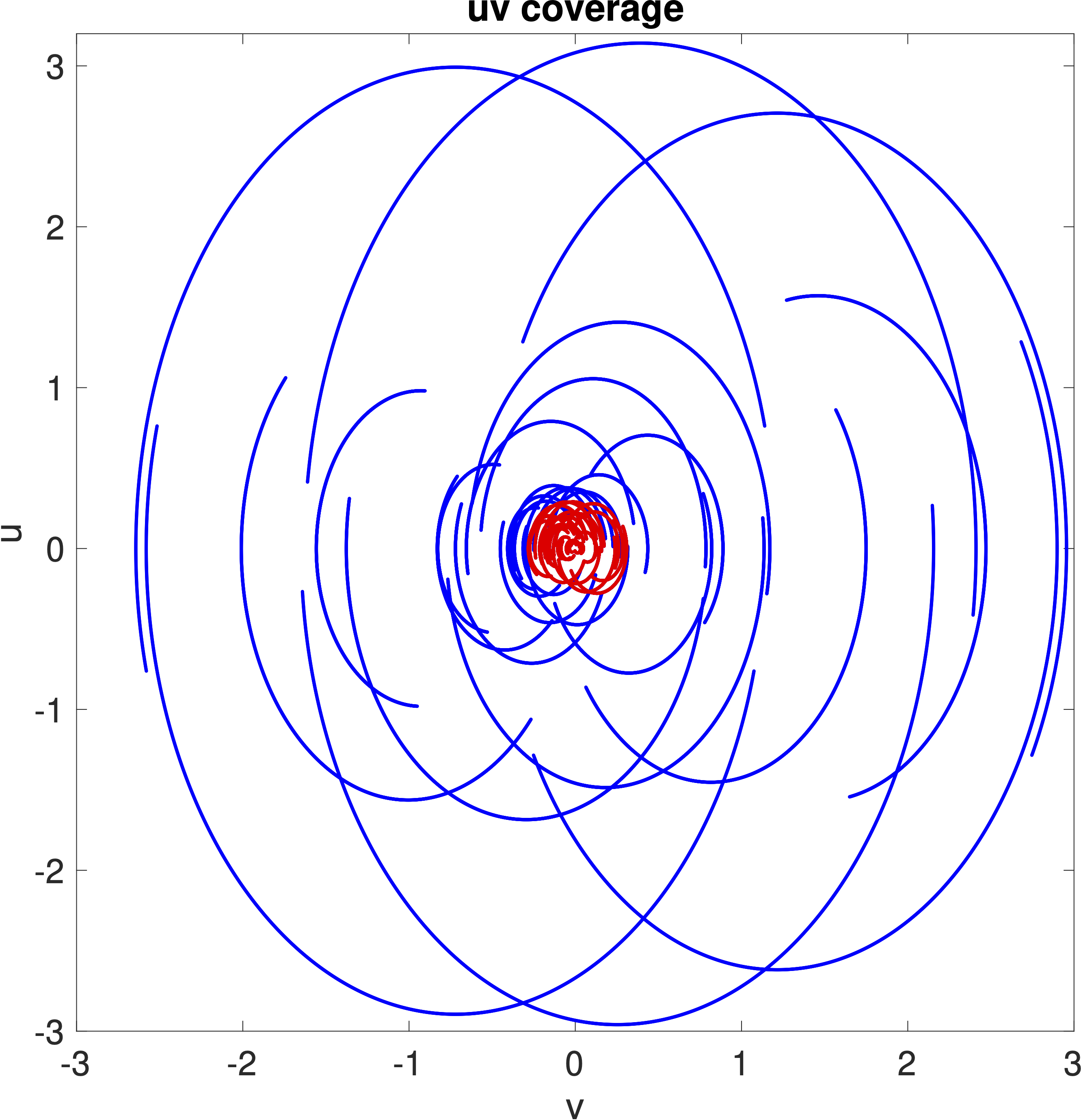}\vspace{-0.5em}
    \caption{$u-v$ coverage for the fine (in \textcolor{blue}{blue}) and coarse level (in \textcolor{red}{red}) when using the MeerKAT telescope \cite{MeerKAT}.\vspace{-2em}}
    \label{fig:uv_coverage}
\end{figure}

\smallskip\noindent \textbf{Construction of $\mathbf{S}$} -- 
To demonstrate the potential of the proposed IML-FISTA for RI imaging, we choose $\mathbf{S}$ in~\eqref{eq:sketching_model} to select low-frequency coefficients in the Fourier $u-v$ coverage, and %
preserving the ellipsis arcs  %
(i.e. corresponding to antenna pairs selecting low-frequency components). %
An example with a $u-v$ coverage simulated from a subset of $64$ antennas of the MeerKAT telescope \cite{MeerKAT} is displayed in Figure~\ref{fig:uv_coverage}, where $\mathbf{S}$ selects $m/2$ coefficients (in red) out of the $m=10,080,000$ total observations. %
\begin{figure*}[ht]
    \centering 
     \includegraphics[trim={0em 0em 0em 0em},clip,width=0.28\textwidth]{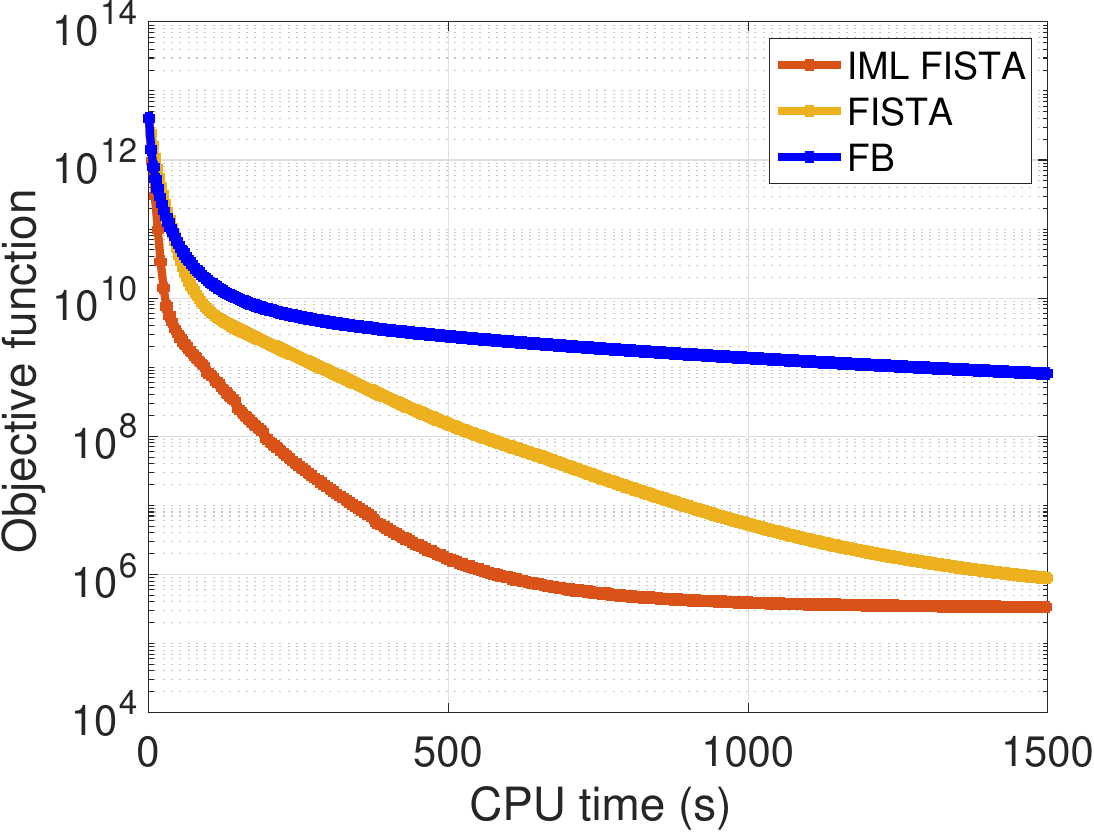} \hspace{1em} \includegraphics[trim={0em 0em 0em 0em},clip,width=0.27\textwidth]{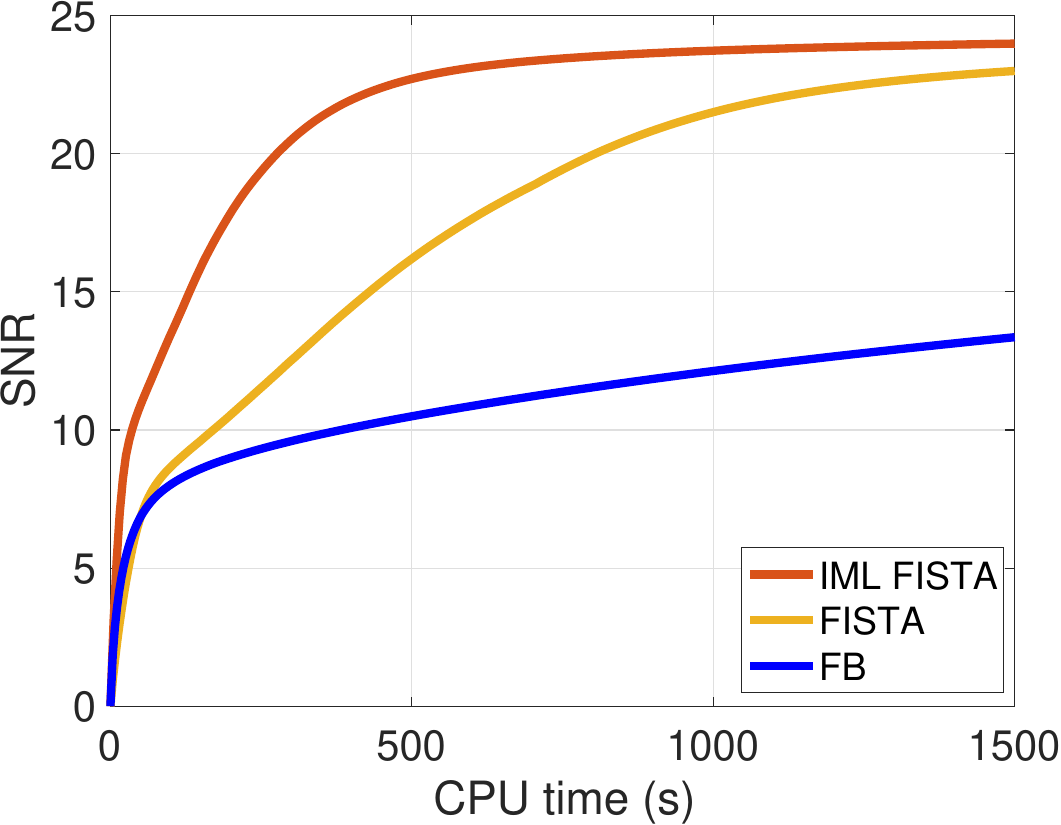}
    \hspace{1em}
    \includegraphics[width=0.27\textwidth]{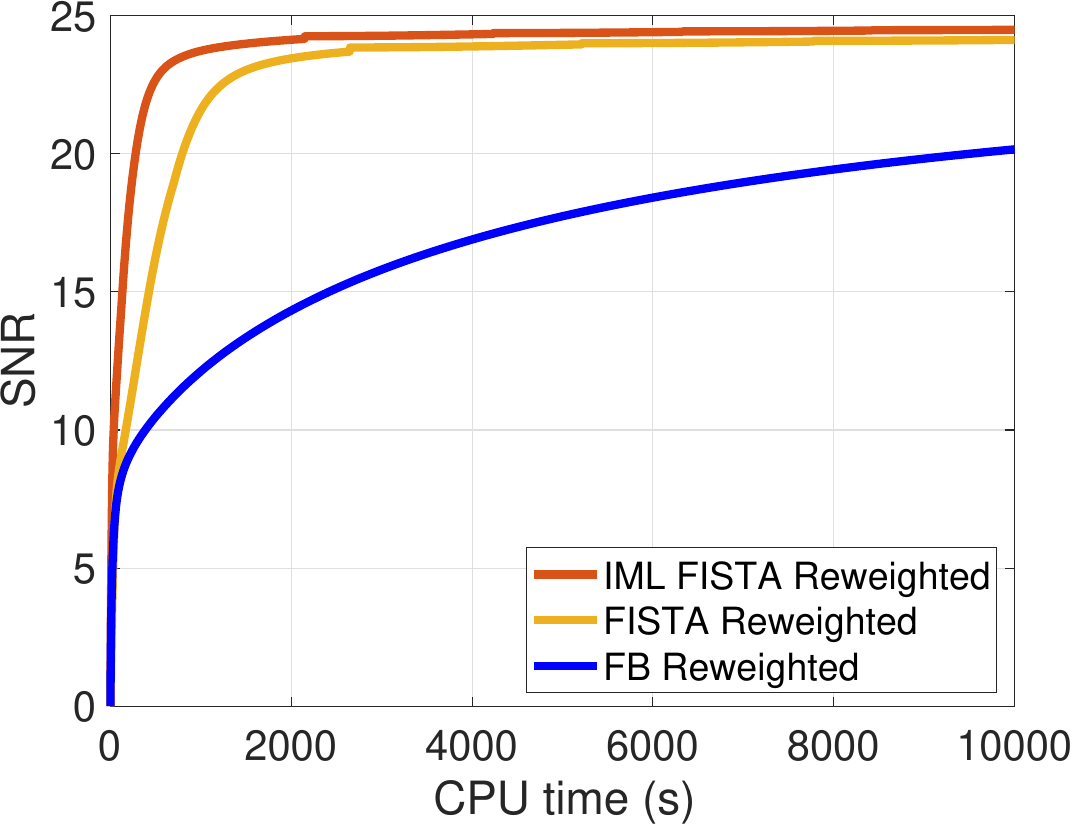}
    \vspace{-0.5em}
    \caption{Evolution  of the objective function values \textbf{(left)}, of the SNR \textbf{(middle)} of the iterates  produced by FB algorithm, FISTA and IML-FISTA with respect to the CPU time when solving Problem \eqref{min:fine-gen} (each algorithm had a CPU time budget of 1500 seconds).  Evolution of the SNR for the three algorithms when we involve the complete reweighting procedure \textbf{(right)} (CPU time budget of 10000 seconds).\vspace{-0.5em}}\label{fig:M31_single_round_FSNR}
\end{figure*}
In other words we keep the visibilities produced by pair of antennas with the smallest distance to each others in the physical world. %
This choice is also based on the fact that most of the signal energy is usually concentrated around low frequencies \cite{carrillo2012sparsity} to accelerate the reconstruction of the image, a common idea in RI imaging \cite{onose2017accelerated}. %

Formally $\mathbf{S}$ is 
a sub-sampling operator that selects a subset $J \subseteq \{1,\ldots,m\}$ of the available visibilities.  %
We then construct %
$\mathbf{\Phi}_H = \mathbf{S\Phi}$ to map the DFT of the image $x$ to $\mathbf{S} y$. %
\begin{figure*}[ht]
    \begin{center}
        \setlength{\tabcolsep}{8pt}
        \setlength{\extrarowheight}{-10pt}
        \begin{tabular}{ccccccc}
        & \ftn$\mathbf{-0.22}$ \textbf{dB} - $277$ s & \ftn$\mathbf{0.01}$ \textbf{dB} - $579$ s & \ftn$\mathbf{0.40}$ \textbf{dB} - $885$ s& \ftn$\mathbf{0.80}$ \textbf{dB} - $1195$ s& \ftn$\mathbf{1.18}$ \textbf{dB} - $1495$ s & \\
        \begin{turn}{+90} \quad \quad FB \end{turn} & \includegraphics[trim={0 0 0 0},clip,width=0.11\textwidth]{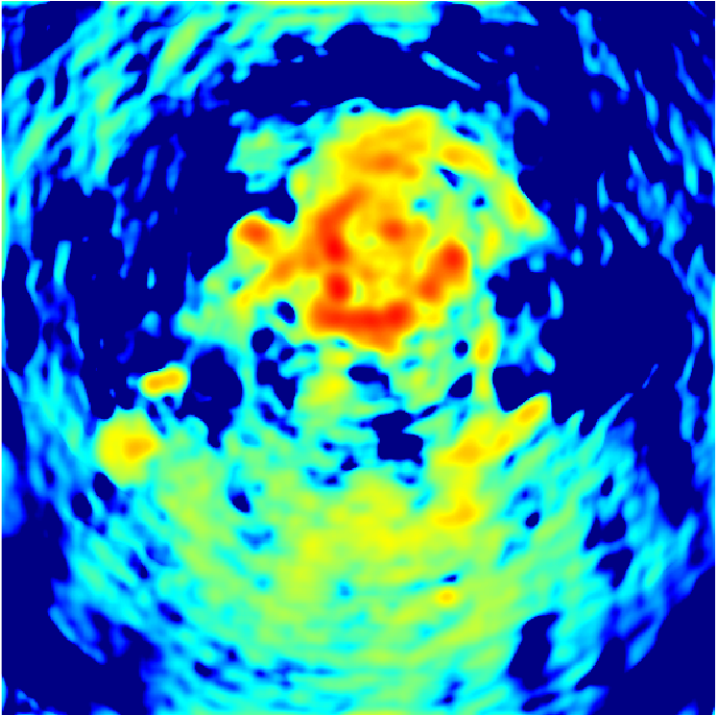} & \includegraphics[trim={0 0 0 0},clip,width=0.11\textwidth]{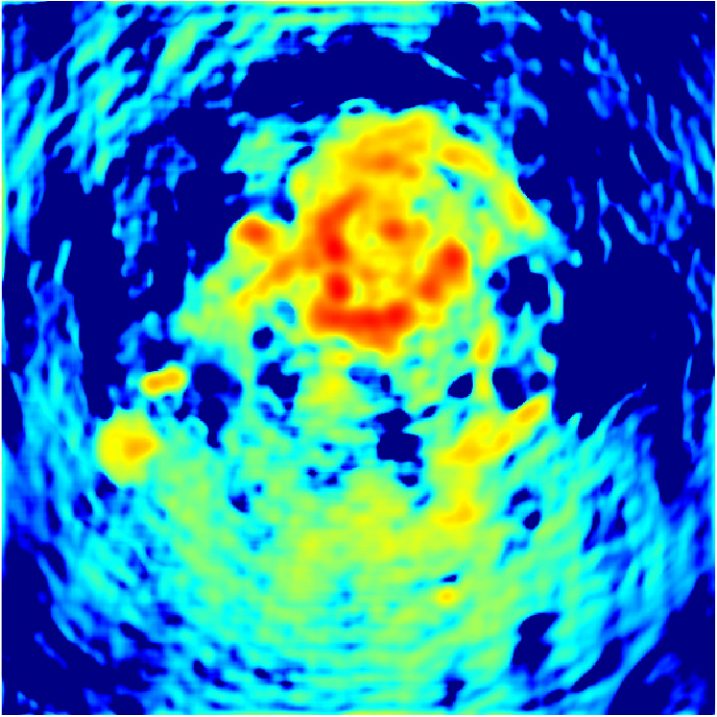} & \includegraphics[trim={0 0 0 0},clip,width=0.11\textwidth]{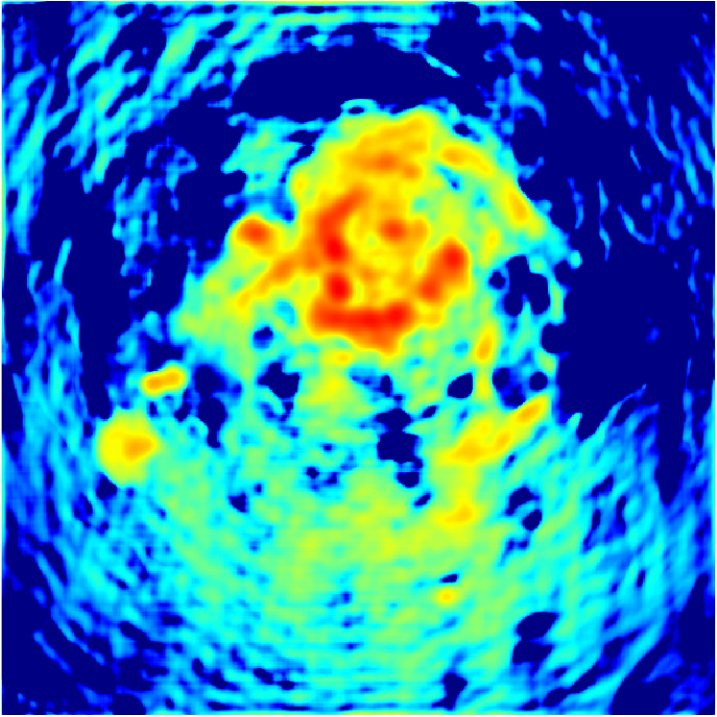} & \includegraphics[trim={0 0 0 0},clip,width=0.11\textwidth]{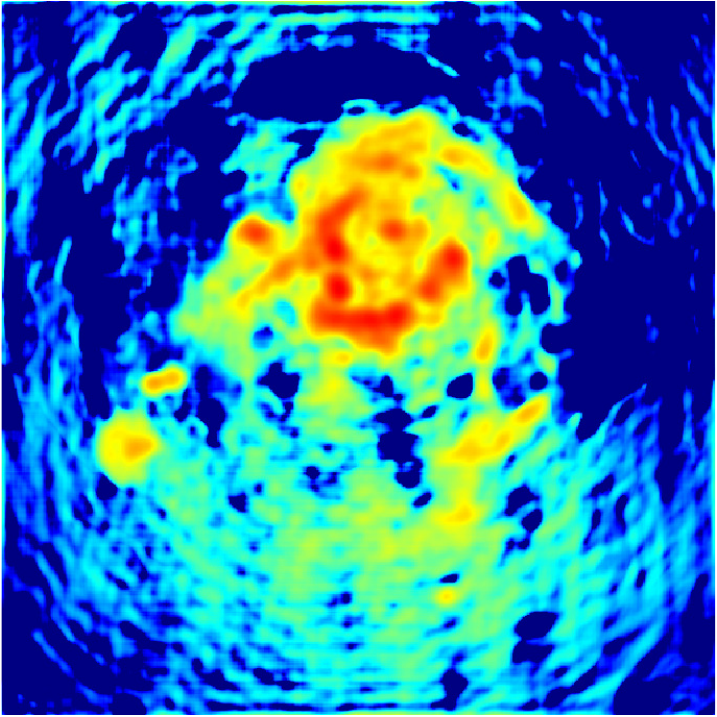} 
         & \includegraphics[trim={0 0 0 0},clip,width=0.11\textwidth]{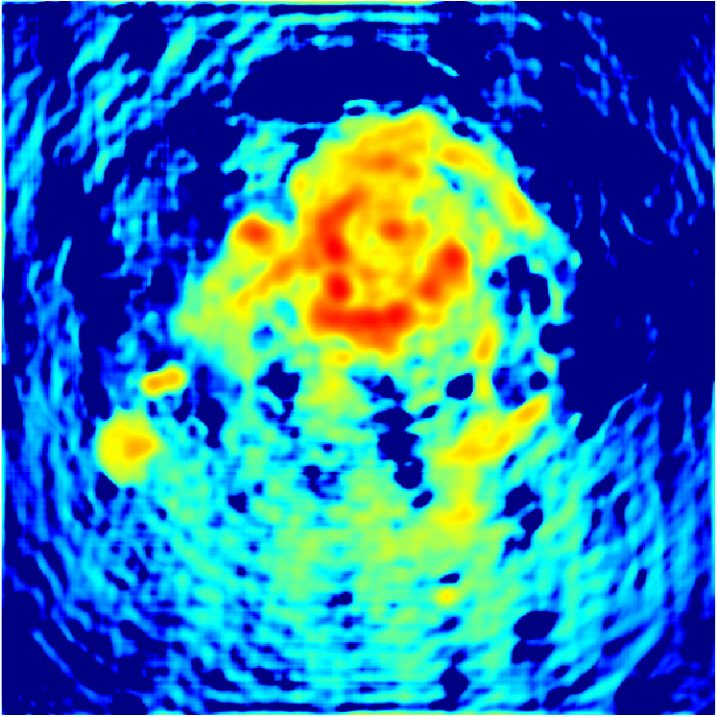} & \multirow{2}{*}{\includegraphics[trim={0 0 0 0},clip,height=0.25\textwidth,width=0.04\textwidth]{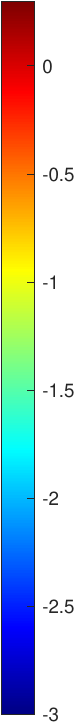}}
         \\
         & \ftn$\mathbf{0.98}$ \textbf{dB} - $303$ s& \ftn$\mathbf{3.49}$ \textbf{dB} - $575$ s& \ftn$\mathbf{8.29}$ \textbf{dB} - $897$ s& \ftn$\mathbf{11.26}$ \textbf{dB} - $1187$ s& \ftn$\mathbf{13.26}$ \textbf{dB} - $1474$ s &
         \\
         \begin{turn}{+90} ~~ FISTA \end{turn} & \includegraphics[trim={0 0 0 0},clip,width=0.11\textwidth]{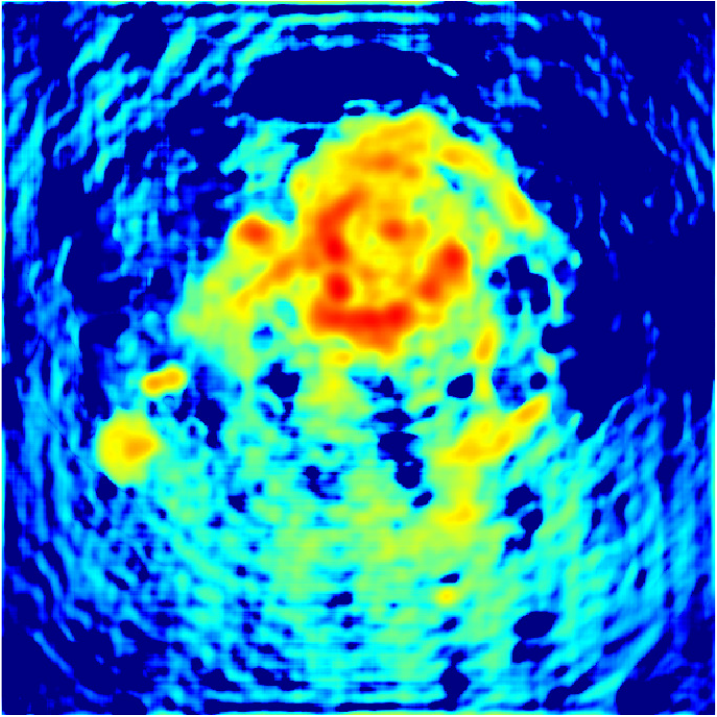} & \includegraphics[trim={0 0 0 0},clip,width=0.11\textwidth]{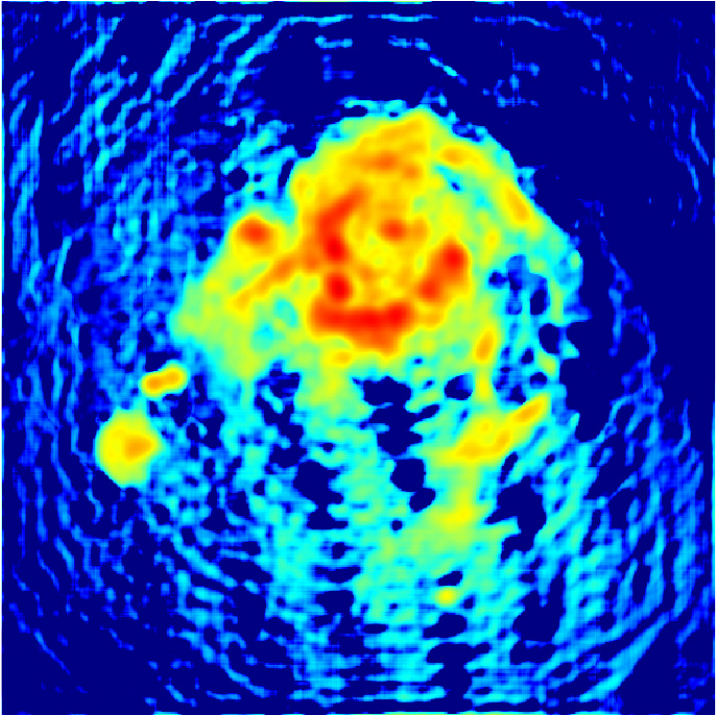} & \includegraphics[trim={0 0 0 0},clip,width=0.11\textwidth]{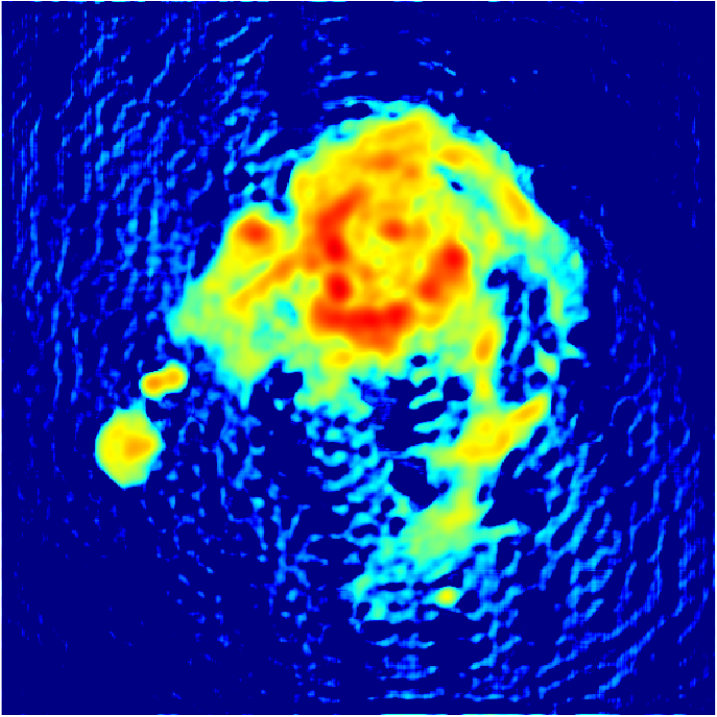} & \includegraphics[trim={0 0 0 0},clip,width=0.11\textwidth]{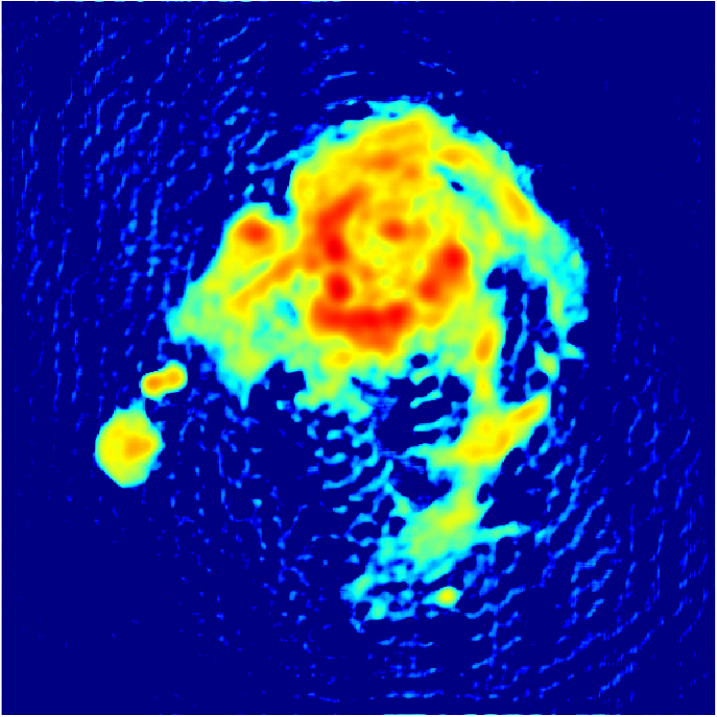} &  \includegraphics[trim={0 0 0 0},clip,width=0.11\textwidth]{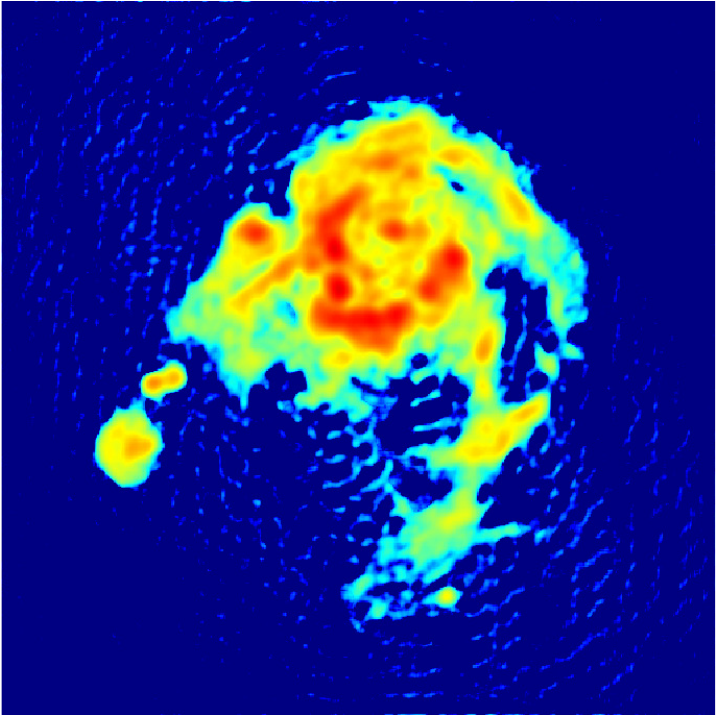}  & %
         \\
         & \ftn$\mathbf{7.63}$ \textbf{dB} - $300$ s & \ftn$\mathbf{11.57}$ \textbf{dB} - $458$ s & \ftn$\mathbf{15.45}$ \textbf{dB} - $769$ s & \ftn$\mathbf{16.89}$ \textbf{dB} - $1060$ s& \ftn$\mathbf{17.83}$ \textbf{dB} - $1497$ s &
         \\
         \begin{turn}{+90} ~IML-FISTA \end{turn} & \includegraphics[trim={0 0 0 0},clip,width=0.11\textwidth]{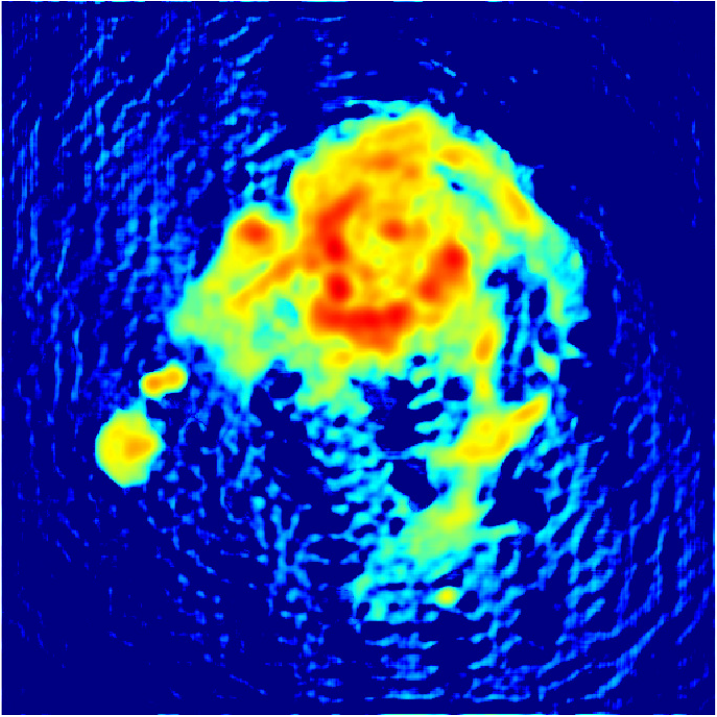} & \includegraphics[trim={0 0 0 0},clip,width=0.11\textwidth]{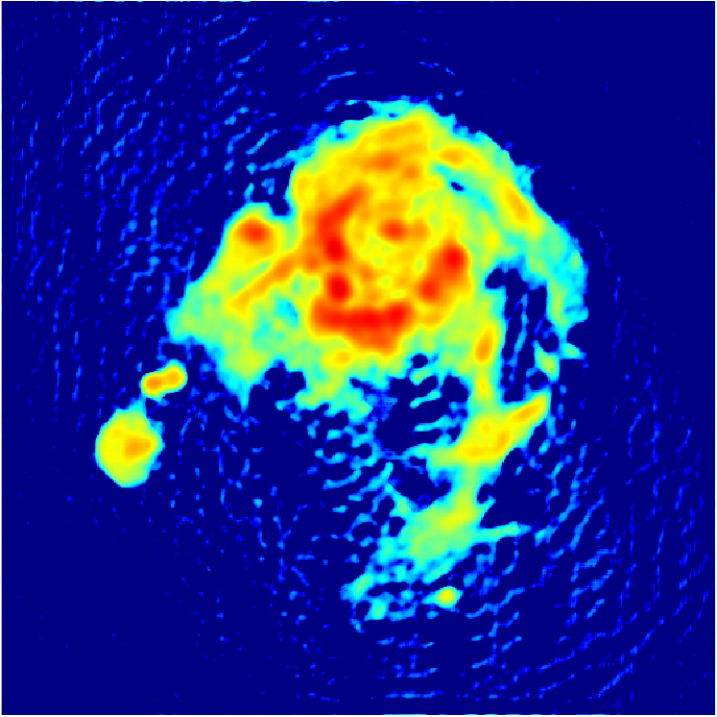} & \includegraphics[trim={0 0 0 0},clip,width=0.11\textwidth]{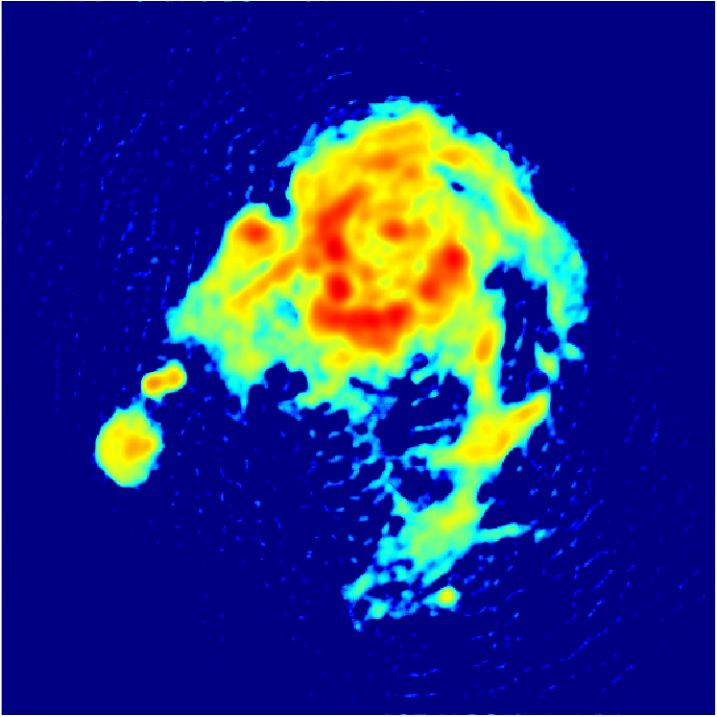} & \includegraphics[trim={0 0 0 0},clip,width=0.11\textwidth]{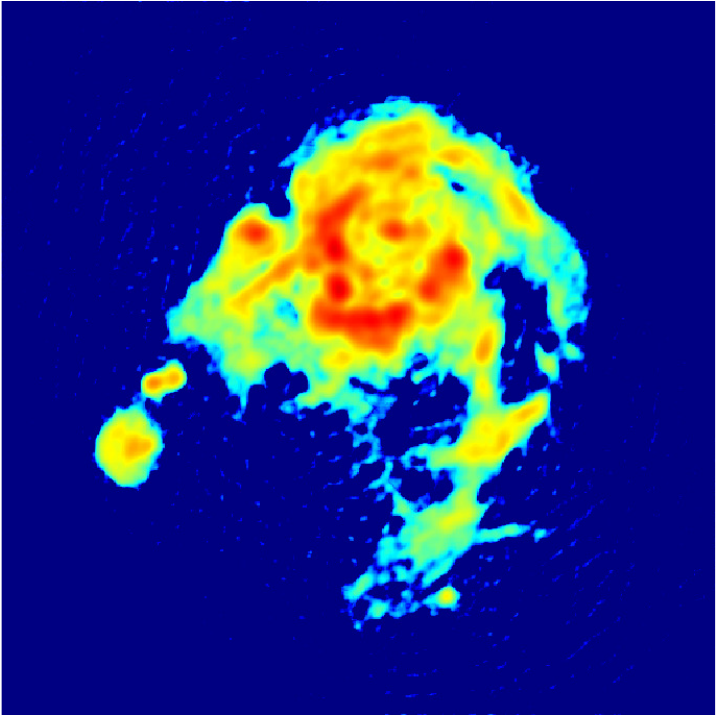} 
         & \includegraphics[trim={0 0 0 0},clip,width=0.11\textwidth]{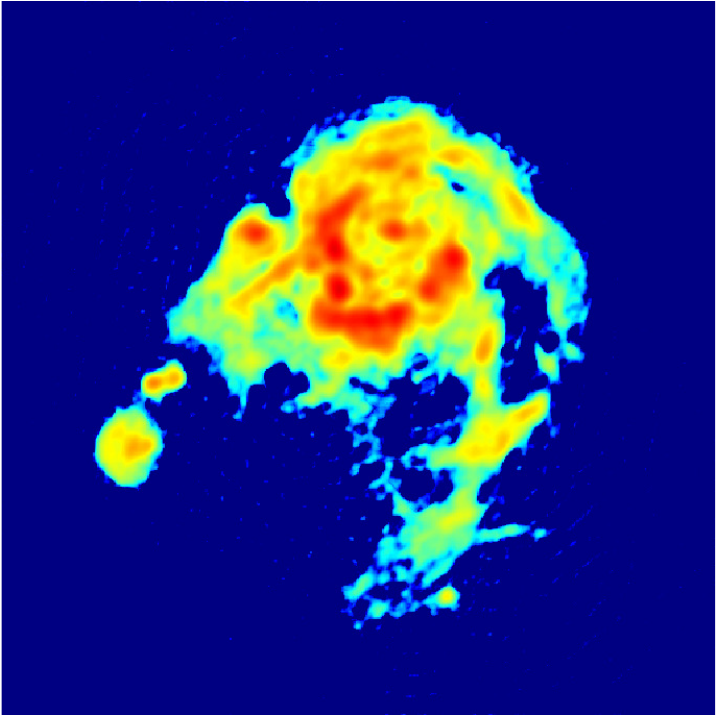} & %
        \end{tabular}
    \end{center} \vspace{-1em}
    \caption{Reconstruction in log scale of a region of the M31 galaxy by FB (top row) FISTA (middle row) and IML-FISTA (bottom row) at equivalent CPU times. The legend on top of each thumbnail reads as follows: \textbf{log SNR} in dB - CPU time in seconds. log SNR = SNR($\log_{10}(10^3 x + 1)/3, \log_{10}(10^3 x_{\text{truth}} + 1)/3$).\vspace{-1.5em}}
    \label{fig:M31_thumbn_1500}
\end{figure*}

\smallskip\noindent\textbf{Choice of coarse model} --
We choose $R_{H,\gamma} = 0$, i.e. the coarse level is not regularized explicitly. This choice is due to the fact that through $\mathbf{S}$ we are only working with low-frequencies, and we observed that adding a coarse regularization in this case was not making a quantitative difference in practice. Thus, the coarse objective function~\eqref{eq:F_H} boils down to
\begin{equation}
    \label{eq:coarse_level}
   (\forall x \in \RR^n) \quad F_H(x) :=\frac{1}{2}\Vert \mathbf{S}\mathbf{\Phi} x-\mathbf{S}y \Vert^2 + \langle v_{H},x \rangle.
\end{equation}
Hence, %
the coarse model is still guided by the fine level regularization through $v_{H,k}$. 

Regarding the smoothing of $R(\cdot, \mathbf{W}_i)$ in~\eqref{eq:v_Hk_inexact1}%
, we choose to only smooth the SARA weighted-$\ell_1$ regularization %
without enforcing the first order coherence with respect to $\iota_{\RR^n_+}(\cdot)$. %
In practice we have not observed unfeasible coarse iterates.

\smallskip\noindent \textbf{Multilevel setup} -- For all of the performed experiments, our ML algorithm will consist of 3 levels. Each coarse level has half of the available measurements of its corresponding fine level. At each level we perform $p=5$ iterations of gradient descent in algorithm~\ref{alg:IMLFISTA} with $F_H$ given in \eqref{eq:coarse_level}. %

\section{Numerical experiments}
\label{sec:numerical_experiments}

\subsection{Simulated data}
We use a subset of $64$ antennas from the MeerKAT array \cite{MeerKAT}. Each antenna pair acquires $5,000$ visibilities, leading to a total of $m = 10,080,000$ observations (see Figure~\ref{fig:uv_coverage} for the resulting Fourier coverage).  
In our simulations, we use a simulated image of the M31 galaxy\footnote{Image available \href{https://casaguides.nrao.edu/index.php?title=Sim_Inputs}{here}.} of dimension $n= 512 \times 512$.
The measurements 
are obtained as per equation~\eqref{pb:inv}, where $\epsilon \in \mathbb{C}^m$ is a realization of a centered white Gaussian noise with variance $\sigma = 0.007$, so that the input Signal-to-Noise-Ratio (SNR) is equal to 19~dB in the visibility domain. %

\subsection{
Minimization comparison without reweighting %
}
In this section we will compare three optimization methods for solving Equation \eqref{min:fine-gen}%
: FB, FISTA, and IML-FISTA. 
Each algorithm is given a budget of CPU time to reach the best reconstruction ($\lambda$ is chosen via grid search). Our main goal is to demonstrate that IML-FISTA is faster than FISTA %
to solve this problem. First and foremost we are interested in the quality of the reconstruction so we will plot two criteria to validate the performances of our algorithm: 
the objective function and the SNR evolution with respect to the CPU time, in Figures~\ref{fig:M31_single_round_FSNR} left and middle, respectively.
As one can see IML-FISTA outperforms both FISTA and FB algorithms for a single round of convex optimization.
We further provide reconstructions obtained with the three methods for visual inspection in Figure \ref{fig:M31_thumbn_1500}, at given CPU computation times $\{\approx 300\text{s}, \approx 500\text{s}, \approx 800\text{s}, \approx 1,100\text{s}, \approx 1,500\text{s}\}$. %
\subsection{%
Minimization comparison for uSARA %
}
We now focus on solving the complete uSARA problem. As we solve a sequence of optimization problems \eqref{eq:reweight_step} that will be different for each optimization method, the easiest way to evaluate each method is to %
only compare FB, FISTA and IML-FISTA on the SNR evolution with respect to the CPU time computation. %
The results are shown in Figure~\ref{fig:M31_single_round_FSNR}-right.
One can see %
that at each "reweighting" step a small jump in the SNR of the iterates occurs (for both FISTA and IML-FISTA) due to the reset of the inertia parameters (for convergence reasons \cite{aujol2015}). With the given coverage we can only slightly improve the SNR of the reconstruction, but nevertheless IML-FISTA reaches an upper bound faster than FISTA.

\section{Conclusion and perspectives}
\noindent \textbf{Conclusion} -- 
We proposed a ML approach for solving the uSARA problem in RI imaging, where the coarse level enables working with low-dimensional data, while the fine model ensures consistency with the full data and promotes averaging sparsity. We have also integrated the resulting IML-FISTA iterations, reminiscent from the ML approach proposed in~\cite{lauga2023}, within a reweighting framework, further enhancing sparsity. %
We have shown through simulations on RI imaging that the proposed IML-FISTA leads to  
impressive acceleration with respect to FB to solve uSARA by exploiting approximations in the observation space of the problem.  
Our method shows promising results on simulations when integrated within the reweighting procedure. %

\noindent \textbf{Perspectives} -- 
The proposed ML acceleration for uSARA being very promising, we have identified future research directions to better assess its potential for RI imaging.

On the one hand, there exist approaches to efficiently handle high-dimensional data based on a parallel implementation of the measurement operator $\mathbf{\Phi}$ \cite{onose2016scalable}. Our method could be coupled with such a parallelisation strategy to benefit at the same time from the dimensionality reduction at the coarse levels and from an efficient parallel implementation of $\mathbf{\Phi}$ at the fine levels. 
Also, preconditioning strategies enabling natural weighting (leveraging the local density of the Fourier sampling) could be considered for comparison and/or further acceleration of the proposed method  \cite{onose2017accelerated}. 

Moreover, sophisticated coarse models for RI could be investigated  to potentially improve the results (for instance sketching approaches \cite{vijay2017fourier}).

Furthermore, connections of ML approaches with CLEAN algorithm \cite{hogbom1974aperture} and its learned version R2D2 \cite{Aghabiglou2023} could be studied. Both methods are built on major-minor cycles reminiscent of matching pursuit. During the minor cycles, an approximate data term is used, ultimately enabling a much smaller number of major cycles (requiring passing through the full data). This is akin to the proposed ML method.%

On the other hand, a few theoretical research directions could be pursued. Leveraging approximation theory as in~\cite{repetti2021variable}, the global convergence of the ML strategy within a reweighting framework could be studied. The ML framework could also be extended to primal-dual algorithms, to enable solving the constrained formulation of SARA, which in this RI context yields better reconstruction quality. \vspace{-0.7em}%

\bibliographystyle{IEEEbib}
\bibliography{strings}

\end{document}